\documentclass{article}

\usepackage{amsmath,amssymb,amscd,theorem,pb-diagram}

{\theorembodyfont{\slshape}\newtheorem{theorem}{Theorem}}
{\theorembodyfont{\slshape}\newtheorem{proposition}[theorem]{Proposition}}
{\theorembodyfont{\slshape}\newtheorem{lemma}[theorem]{Lemma}}
{\theorembodyfont{\slshape}\newtheorem{definition}[theorem]{Definition}}
{\theorembodyfont{\rmfamily}\newtheorem{note}[theorem]{Note}}
{\theorembodyfont{\slshape}}
{\theorembodyfont{\slshape}}
{\theorembodyfont{\rmfamily}}

\newcommand{\Fq}{\ensuremath{\mathbb{F}_q}}

\renewcommand{\O}{\ensuremath{\mathcal{O}}}
\newcommand{\Ot}{\ensuremath{\widetilde{\mathcal{O}}}}
\newcommand{\G}{\ensuremath{\Gamma}}
\newcommand{\g}{\ensuremath{\gamma}}
\newcommand{\og}{\ensuremath{\bar{\gamma}}}
\newcommand{\Q}{\ensuremath{\mathbb{Q}}}

\newcommand{\Qq}{\ensuremath{\mathbb{Q}_q}}

\newcommand{\Zq}{\ensuremath{\mathbb{Z}_q}}
\newcommand{\Qt}{\ensuremath{\mathbb{Q}_2}}
\newcommand{\Zt}{\ensuremath{\mathbb{Z}_2}}

\newcommand{\C}{\ensuremath{\mathbb{C}}}
\newcommand{\Z}{\ensuremath{\mathbb{Z}}}
\newcommand{\N}{\ensuremath{\mathbb{N}}}
\newcommand{\F}{\ensuremath{\mathbb{F}}}
\renewcommand{\S}{\ensuremath{\mathbb{S}}}

\renewcommand{\k}{\ensuremath{\kappa}}

\newcommand{\ord}{\ensuremath{\text{ord}}}

\newcommand{\Res}{\ensuremath{\text{Res}}}
\newcommand{\B}{\ensuremath{\mathcal{B}}}
\newcommand{\ep}{\hfill $\blacksquare$\vspace{\baselineskip}} 

\begin{document}

\title{Point counting in families of hyperelliptic curves in characteristic 2}

\author{ Hendrik Hubrechts\\
\small{Research Assistant of the Research Foundation - Flanders (FWO - Vlaanderen)}\\
\small{Department of mathematics, Katholieke Universiteit Leuven}\\
\small{Celestijnenlaan 200B, 3001 Leuven (Belgium)}\\
\small{\texttt{Hendrik.Hubrechts@wis.kuleuven.be}}}

\maketitle

\begin{abstract}
Let $\bar E_{\G}$ be a family of hyperelliptic curves over $\F_2^{\text{alg cl}}$ with general Weierstrass equation given over a very small field $\F$. We describe in this paper an algorithm to compute the zeta function of $\bar{E}_{\og}$ for $\og$ in a degree $n$ extension field of $\F$, which has as time complexity $\Ot(n^3)$ and memory requirements $\O(n^2)$. With a slightly different algorithm we can get time $\O(n^{2.667})$ and memory $\O(n^{2.5})$, and the computation of $\O(n)$ curves of the family can be done in time and space $\Ot(n^3)$. All these algorithms are polynomial in the genus.
\end{abstract}


\section{Introduction and results}\label{sec:intro}
The problem of counting rational points on curves over finite fields has received much attention during the last decade, and many algorithms have been proposed. For an overview of these results and their relevance we refer to \cite{KedlayaComputingZetaFunctions, VercauterenThesis}. (Hyper)elliptic curves over finite fields of characteristic 2 are particularly interesting due to the fact that computers can work very efficiently with them.

For elliptic curves Mestre has presented an algorithm using the arithmetic geometric mean (AGM) that works in time $\Ot(n^3)$, and Lercier and Lubicz \cite{LercierLubicz} extended and improved it to $\Ot(n^2)$ and very small genus $>1$. Kedlaya presented in \cite{KedlayaCountingPoints} an algorithm to compute the zeta function of hyperelliptic curves of genus $g$ in odd characteristic in time $\Ot(g^4n^3)$ using Monsky-Washnitzer cohomology, and Denef and Vercauteren \cite{DenefVercauteren} extended this to characteristic 2.

On the other hand Lauder \cite{LauderDeformation} and Tsuzuki \cite{TsuzukiKloosterman} introduced deformation in the story of point counting, and in \cite{HubrechtsHECOdd} we followed a suggestion of Lauder to combine deformation with Kedlaya's approach. In this paper we extend this result to characteristic 2, thereby reconciling Denef and Vercauteren's work with deformation.

In \cite{GerkmannHypersurfaces} Gerkmann also considered a deformation approach for elliptic curves in odd characteristic, and due to the easy form of the Weierstrass equation of an elliptic curve in characteristic 2 he was able to add this situation without much effort. For higher genus this equation is however more involved, and as a consequence the theory is technically rather different from the odd characteristic case, although the `big picture' has a similar \emph{esprit}.

We will now present the results proven in this paper. Let $\Fq$ be a finite field with $q=2^a$ elements, $\og\in\F_{q^n}$ for some integer $n$, and $g\geq 1$ an integer. Suppose $\bar{f},\bar{h}\in\Fq[X,\G]$ are in the form described in section \ref{ssec:introDef}, which implies especially that for most $\og$ we get a hyperelliptic curve of genus $g$ over $\F_{q^n}$ of the form
\[\bar E_{\og}: Y^2 +\bar{h}(X,\og)Y=\bar{f}(X,\og).\]
Define $\k:=\max\{\deg_{\G}f,\deg_{\G} h^2\}$. As is mentioned in \cite{DenefVercauteren}, in this matter we have an `average case' and a `worst case'. This means that almost all curves belong to the first case, and some unlucky ones do not. In this paper we will often use the Soft-Oh notation $\Ot$ as defined in \cite{ModernCompAlg}, which is essentially a Big-Oh notation that ignores logarithmic factors.

The main result is the following theorem, to be proven in section \ref{sec:complexity}.
\begin{theorem}\label{thm:princThm}
We can compute deterministically the zeta function of (the projective completion of) $\bar E_{\og}$ using $\Ot(g^{6,376}a^3\k^2n^2+g^{3,376}a^3n^3)$ bit operations and $\Ot(g^5a^3\k n^2)$ bits of memory `on average'. For the worst case scenario one factor $g$ is to be added to the terms with $n^2$ in them.\end{theorem}
By using some faster substitution algorithm it is possible to gain time, at the cost of an increase in memory usage. The result is the following.
\begin{theorem}\label{thm:subcubic}
There exists a deterministic algorithm that computes the zeta function of $\bar E_{\g}$ in $\Ot(g^{6,376}a^3\k^2n^2+g^{3,376}a^2n^{2,667})$ bit operations `on average'. It requires then $\Ot(g^5a^3\k n^2+g^3a^2n^{2,5})$ bits of memory. In the worst case again one factor $g$ has to be added to both first terms.\end{theorem}
Theorem \ref{thm:subcubic} together with the following result and an algorithm quadratic in $n$ for a special situation with a Gaussian normal basis is proven in section \ref{sec:conclusion}. In this theorem we did not pay attention to the dependency of parameters different from $n$.
\begin{theorem}\label{thm:lotscurves}
Given $\O(n)$ parameters $\og_1,\ldots, \og_k\in\F_{q^n}$, it is possible to find the zeta functions of all $\bar E_{\og_i}$ with $\Ot(n^3)$ as time and space requirements.\end{theorem}
The bottom line of this algorithm is that in order to find a curve with some special size by trying a lot of curves, we can count on $\Ot(n^2)$ as the time needed for one curve.\\

This paper is organized as follows. In section \ref{sec:analysis} we provide the theory behind the algorithm, in it is explained the required special form of $\bar f$ and $\bar h$, to which we referred earlier. In section \ref{sec:matrices} we have gathered some necessary results about certain $2$-adic matrices and differential equations of them, required for the algorithm. More precisely, some trick is explained to compute the matrix of the connection and a particularly useful form of the differential equation, the convergence properties of Frobenius are investigated, and an important result about error control is established. The next section gives the algorithm and proves its correctness, and section \ref{sec:complexity} estimates the complexity, thereby proving theorem 1. Finally the last section mentions the improvements noted above, in particular theorems \ref{thm:subcubic} and \ref{thm:lotscurves}.

\section{Analytic theory}\label{sec:analysis}
In this section we will develop an analytic theory which combines the results from \cite{DenefVercauteren} with a deformation. Before we start let us define some notation used throughout the rest of the paper. Let $a$ be a strictly positive integer, then we denote by $\Fq$ the finite field with $q := 2^a$ elements. Let $\Qt$ be the completion of $\Q$ according to the 2-adic norm, and $\Qq$ is the unique degree $a$ unramified extension of $\Qt$. Denote with $\C_2$ the completion of an algebraic closure of $\Q_2$. The rings of integers of $\Qt$ and $\Qq$ are written $\Zt$ respectively $\Zq$. The lift of the Frobenius automorphism on $\Fq$ is given by $\sigma:\Qq\to\Qq$. We extend $\sigma$ by acting as squaring on each appearing variable unless said otherwise. If $k$ is a field, then we mean by $k^{\text{alg cl}}$ an algebraic closure of $k$. The derivative of some expression $\alpha$ with respect to $X$ will be denoted by $\alpha'$, and on the other hand $\frac{\partial \alpha}{\partial \G}$ is written as $\dot{\alpha}$.

\subsection{Introducing the deformation.}\label{ssec:introDef}

Suppose we are given an equation $Y^2+\tilde{h}(X)\cdot Y=\tilde{f}(X)$ over $\Fq$ which defines a hyperelliptic curve of genus $g$. As pointed out in \cite{DenefVercauteren} it is always possible to find in an efficient way an isomorphic curve over $\Fq$ given by $Y^2+\bar{h}(X)\cdot Y=\bar{f}(X)$ subject to the following conditions. The degree of the monic polynomial $\bar{f}$ is $2g+1$ and $\bar{h}$ is nonzero of degree at most $g$. If we factor $\bar{h}$ in its monic irreducible factors over $\Fq$, $\bar{h}(X)=\bar{c}\prod_{i=1}^s\bar{h}_i^{r_i}(X)$ with $\bar{h}_i$ irreducible, $r_i\neq 0$ and $\bar{c}\in\Fq^\times$, define then $\bar{H}(X) := \prod_{i=1}^s\bar{h}_i(X)$, the product of the irreducible factors of $\bar{h}$. We require now that $\bar{f} = \bar{H}\cdot \bar{Q}_{\bar{f}}$ where $\bar{H}$ and $\bar{Q}_{\bar{f}}$ are relatively prime.

Define $\tilde D:=\max{r_i}$ so that $\bar{h}$ is a divisor of $\bar{H}^{\tilde D}$, and let $\bar{Q}_{\bar{h}}$ be such that $\bar{h}\cdot \bar{Q}_{\bar{h}} = \bar{H}^{\tilde D}$. Now we can lift the $\bar{h}_i$ and $\bar{Q}_{\bar{f}}$ to $h_i$ and $Q_f$ over $\Zq$ such that they remain monic and the projection modulo 2 equals the original polynomials. As a consequence we can also define $H$, $h$ and $Q_h$ with the same properties as in the finite field case.\\

To introduce the deformation parameter $\G$ in the resulting equations, we allow $Q_f$ and the $h_i$ to be polynomials in $\Zq[X,\G]$ such that they remain monic in $X$. Let $r(\G)$ be the resultant of $H$ and $Q_f\cdot\frac{\partial H}{\partial X}=Q_fH'$ with respect to $X$. Then we require $r(\G)$ to be a polynomial for which $r(0)$ does not reduce to zero modulo 2, or equivalently $\G=0$ gives a hyperelliptic curve modulo 2. The resultant determines for which parameters the result is a hyperelliptic curve, therefore we define the following subset of the set $\text{Teich}(\F_2^\text{alg cl})$ of Teichm\"{u}ller lifts of $\F_2^\text{alg cl}$:
\[\S := \left\{\g\in\text{Teich}(\F_2^\text{alg cl})\ \right|\left. \vphantom{\text{Teich}(\F_2^\text{alg cl})}\ r(\g)\not\equiv 0\bmod 2\right\}.\]
\begin{lemma}\label{lem:SGoodParams}
For $\g\in\S$ the projected equation $Y^2+\bar{h}(X,\og)\cdot Y=\bar{f}(X,\og)$ defines a hyperelliptic curve $\bar{E}_{\og}$ over $\F_2^\text{alg cl}$, with an equation of the form mentioned above.
\end{lemma}
\textsc{Proof.} It is enough to show for a Teichm\"{u}ller lift $\g$ that $\bar{E}_{\og}$ has no affine singularities iff $\g\in\S$. Computing the system of partial derivatives yields immediately that the existence of an affine singularity $(\bar{x},\bar{y})$ implies that $\bar{H}(\bar{x})=\bar{h}(\bar{x})=\bar{f}(\bar{x})=\bar y=0$ and $\bar{f}'(\bar{x})=0$, and vice versa: these equalities give an affine singularity. As $\bar{f}'=\bar{Q}_{\bar{f}}'\bar{H}+\bar{Q}_{\bar{f}}\bar{H}'$ we conclude that equivalently the system $\bar{H}=\bar{Q}_{\bar{f}}\bar{H}'=0$ has no solutions, which in turn is equivalent to $\Res_X(\bar{H},\bar{Q}_{\bar{f}}\bar{H}')\neq 0$. The fact that the equation has the right structure is trivially checked.\ep

The constructions above fail when $\bar h$ is a constant, in which case $\bar h\neq 0$ is equivalent with $\bar{E}_{\bar{\g}}$ being hyperelliptic for every $\bar{\g}$. In this situation no resultant is needed, and for example $S$ defined below will simply be $\Qq[\G]^\dagger$. We will not always mention the simplifications needed for this special case. The convention $\tilde D=1$ in this case is best suited for the estimates further on.\\

As final definitions, let $\rho:=\deg_\G r(\G)$, $s := \deg_X(H)$ and $\kappa := \max\{\deg_\G f,$ $\deg_\G h^2\}$ as defined before, and $\eta := \deg_\G H$.

\subsection{The overconvergent structures.}\label{ssec:overcvgtStr}

We define as in \cite{HubrechtsHECOdd} the necessary overconvergent structures. For $r=\sum_{i=0}^\rho r_i\G^i$ let $\rho'$ be the largest index for which $\ord(r_{\rho'})=0$, and define $\tilde r=\sum_{i=0}^{\rho'}r_i\G^i$. Hence $\tilde r\equiv r\mod 2$ and if the leading term of $r$ is a unit in $\Z_q$ we have simply $\tilde r= r$. The ring $S$ will be the equivalent of the field $\Qq$ in Denef and Vercauteren's approach.
\[S := \Qq\left[\G,\frac 1{\tilde r(\G)}\right]^\dagger = \left\{ \sum_{k\in\Z}\frac{b_k(\G)}{\tilde r(\G)^k}\ \right|\ (\forall k)\ b_k(\G)\in\Qq[\G], \]\[\qquad\qquad\qquad\qquad\qquad\qquad\qquad
\left.\deg b_k(\G)<\rho' \text{ and }\vphantom{\sum_{k\in\Z}\frac{b_k(\G)}{r(\G)^k}} \liminf_k\frac{\ord(b_k)}{|k|}>0\right\}.\]
The last inequality in this definition is equivalent with the existence of real constant $\delta>0$ and $\varepsilon$ such that for all $k$ we have $\ord(b_k)\geq \delta\cdot|k|+\varepsilon$.
As proven in \cite{HubrechtsHECOdd} the fact that $\tilde r$ is `monic' implies that a general element of $S$ can be represented as $\sum_{i=0}^\infty a_i\G^i+\sum_{j=1}^\infty\frac{b_j(\G)}{\tilde r(\G)^j}$ where we have $\liminf_i{\ord(a_i)}/{|i|}>0$ and $\liminf_j
{\ord(b_j)}/{|j|}>0$. If $\tilde r$ is a constant then of course $S=\Q_q[\G]^\dagger$, and the parts with denominators disappear everywhere. We will not always mention this special case. The equality \[\frac 1r=\frac 1{\tilde r}\sum_{i=0}^\infty\left( -\frac{r-\tilde r}{\tilde r}\right)^i\] combined with the fact that $r-\tilde r\equiv 0 \mod 2$ shows that $1/r\in S$. It is worth remarking that $\S$ doesn't change if defined using $\tilde r$, and $S$ can still be interpreted as consisting of the analytic functions defined over $\Q_q$ and convergent in a disk strictly bigger than the unit disk with small disks removed around the Teichm\"{u}ller lifts not in $\S$. The following important lemma is also proven in \cite{HubrechtsHECOdd} and gives us control over the substitution of some $\g\in\S$ in an element of $s\in S$. Remark that $s(\g)$ always converges.
\begin{lemma}\label{lem:subGinS}
Let $s(\G)=\sum_{k\in\Z}b_k(\G)/\tilde r(\G)^k\in S$. Suppose we have for infinitely many $\g\in\S$ that $\ord(s(\g))\geq\alpha$ for some real number $\alpha$, then also for every $k\in\Z$ we get $\ord(b_k)\geq \alpha$.
\end{lemma}
Now we can define what will be the analogue of the dagger ring $A^\dagger$. The last condition may look quite terrifying, but is a technical condition that implies that the sum $\sum_ks_{ik}$ is convergent and again an element of $S$.
\[T := \frac{\Qq\left[ \G,\frac 1{r(\G)},X,Y,\frac 1{H(X,\G)} \right]^\dagger}{(Y^2-hY-f)}= \left\{ \sum_{k\in\Z}\frac{\sum_{i=0}^{s-1}s_{ik}X^i+\sum_{i=0}^{s-1}s_{ik}'X^iY} {H(X,\G)^k}\ \right|\]\[
(\forall i,k)\ s_{ik}^{(')}\in S,\ \ (\forall i)\ \exists C\in\Qq, \delta>0 \text{ s.t. with } s_{ik}^{(')}=\sum_{j\in\Z}\frac{s_{ikj}^{(')}(\G)}{\tilde r^j}\text{ we have }\]\[\left.(\forall k,j)\ \ord(C \cdot s_{ikj}^{(')})\geq \delta\cdot(|k|+|j|)\vphantom{\sum_{k\in\Z}\frac{\sum_{i=0}^{s-1} s_{ik}X^i+\sum_{i=0}^{s-1}s_{ik}'X^iY} {H(X,\G)^k}}\right\}.\]
In the case that $H$ is a constant the sum over $k$ is restricted to $k\in\Z_{\leq 0}$, which means simply that in this case no denominators with respect to $X$ occur in a general element of $T$. We will write such a general element of $T$ as
\[\sum_{k\in\Z}\frac{U_k(X,\G)+Y\cdot V_k(X,\G)}{H^k},\]
where $\deg_XU_k,V_k\leq s-1$ and $\liminf_k\frac{\ord(U_k,V_k)}{|k|}>0$. It is not hard to expand the proof of lemma 14 given in \cite{HubrechtsHECOdd} such that it yields that $T$ is an $S$-algebra.

Let $\g\in\S$ with $\og\in\F_{q'}$ such that $q'$ is minimal and $\F_q\subset\F_{q'}$. Then we can substitute $\g$ for $\G$ in the above construction of $T$ resulting in the vector space $T(\g)$ over $\Q_{q'}$. We have just as in the odd characteristic case that $T(\g)=A^\dagger$ with $A^\dagger$ as defined in \cite{DenefVercauteren} for the curve $Y^2-h(X,\g)Y-f(X,\g)=0$.\\

We define the derivative with respect to $X$ on $T$ by interpreting $Y$ in terms of $X$. Using the equation in its original form and as $(2Y+h)^2=4f+h^2$ this yields
\begin{equation}\label{eq:derivativeY}Y'=\frac{f'-h'Y}{2Y+h}\cdot\frac{2Y+h}{2Y+h} = \frac{f'h-2fh'+(2f'+hh')Y}{4f+h^2}.\end{equation}
We indeed have that $Y'\in T$ and can hence define the differential $d := T\to TdX: t\mapsto \frac{\partial t}{\partial X}dX$. Let $\imath$ be the hyperelliptic involution $X\mapsto X$ and $Y\mapsto -Y-h(X,\G)$ on $T$, then we have the following central proposition.
\begin{proposition}\label{thm:FreeQuotientModule}
The module $H_{MW} := \frac{TdX}{dT}$ splits into two eigenspaces under $\imath$, namely $H_{MW}^+$ for eigenvalue $+1$ and $H_{MW}^-$ for $-1$. Both are free $S$-modules with basis respectively $\{\frac{X^i}{H}dX\}_{i=0}^{s-1}$ and $\B := \{X^iYdX\}_{i=0}^{2g-1}$.\end{proposition}
If $H$ is a constant, the first basis is empty, or equivalently $H_{MW}^+$ is trivial.\\
\textsc{Proof.} Let $(U+VY)H^{-k}$ be a general term of an element of $T$. Writing $U+VY=\tilde{U}+\tilde{V}(Y+h/2)$ and computing $\imath(Y')=-Y'-h'$ we can readily check that $\imath\circ d=d\circ\imath$, which gives the isomorphism $H_{MW}\cong H_{MW}^+\oplus H_{MW}^-$. Here $\tilde{U}$ gives the first part and $\tilde{V}(Y+h/2)$ the second part.
The linear independence of the elements of the bases can be proven with lemma \ref{lem:subGinS}. Indeed, suppose we have a linear relation $\sum s_ib_i=0$ for basis elements $b_i$ and $s_i\in S$ where $s_j\neq 0$. The lemma then implies the existence of some $\g\in\S$ such that $s_j(\g)\neq 0$, which gives a nontrivial relation $\sum s_i(\g)b_i=0$ in the case without deformation, in contradiction with \cite{DenefVercauteren}.\\

In order to reduce a general element $$\sum_{i\in\Z}U_i(X,\G)dX/H^i+\sum_{j\in\Z}V_j(X,\G)YdX/H^j$$ of $T$, we consider as in \cite{DenefVercauteren} four cases. First, the part with $i\leq 0$ is an exact form, as integrating does not change the overconvergence property. Second, for $i>0$ we have the following formulae from \cite{DenefVercauteren}, where $r_1(\G) := \Res_X(H,H')$, a divisor of $r(\G)$. Write $x^k r_1(\G)=A(X,\G)H+B(X,\G)H'$, and by computing the differential $d(B/H^{i-1})$ we find
\[\frac{x^k}{H^i}\,dX\equiv \frac{1}{r_1}\left( \frac{A}{H^{i-1}}-\frac{B'}{(i-1)H^{i-1}}\right)dX.\]
Repeating this we end with $i=1$ --- which cannot be reduced further, ergo the first basis of the proposition --- and an expression without denominators $H$ which is an exact form. Next, for the part with $j\leq 0$ we can use the following congruence
\begin{equation}\label{eq:congruence}
\left(X^j(2f'+hh')+\frac j3X^{j-1}(4f+h^2)\right)YdX\equiv 0,
\end{equation}
which has leading coefficient $2(2g+1)+4j/3\neq 0$.

Finally we consider the case $j>0$. Let $h=HQ_H$, then by writing $x^k r(\G)=AH+BQ_fH'$ we have
\[\frac{x^k}{H^j}YdX\equiv \]\[\frac 1r\left(\frac{A}{H^{j-1}}+\frac{B(jH'Q_H^2-6Q_f'-3Q_Hh')-B'(4Q_f+Q_Hh)} {(6-4j)H^{j-1}}\right)YdX+\frac{IdX}{rH}.\]
Here the last part $IdX/rH$ is some differential, invariant under the hyperelliptic involution.

Although the above formulae allow us to reduce elements of $T$, they do not guarantee a priori that the reduced elements and the exact differentials appearing are overconvergent. We will prove this for the case $j\leq 0$, the other cases are similar --- the basic idea being that the orders decrease with only logarithmic behavior and $\deg_{\G}$ and `$\deg_r$' increase at most linearly. Let our element of $T$ be given in the form $\sum_{j\geq 0}s_j(\G)X^jYdX$, where $s_j(\G)=\sum_i s_{ij}(\G)\tilde r(\G)^i$ and --- if necessary after multiplying with some constant --- $\ord_2(s_{ij})\geq \delta(j+|i|)$ for some $\delta>0$. It follows immediately from formula (\ref{eq:congruence}) that if $X^jYdX=\sum_bf_b^j(\G)b+dg$, where $b$ runs over $\B$, we have $\deg_{\G}f_b^j\leq \kappa j$, whereas lemma 2 from \cite{DenefVercauteren} and lemma \ref{lem:subGinS} above give that $\ord_2 f_b^j\geq -\left(3+\log_2(j+g+1)\right)$. It is clear that as the coefficients of the original expression grow linearly, we can ignore this logarithmic surplus of the reductions and hence suppose that the $f_b^j$ are integral. If we write
\[\sum_{j=0}^\infty s_j(\G)X^jYdX\equiv \sum_{b\in\B}f_b(\G)b,\]
then we must show that $\sum_js_jf_b^j\in S$. We will prove that with $s_jf_b^j=\sum_{t}\alpha_{tj}(\G)\tilde r(\G)^t$ an inequality $\ord_2(\alpha_{tj})\geq \varepsilon(|t|+j)$ holds, after which fact 10 and lemma 11 from \cite{HubrechtsHECOdd} give the result. Expanding $f_b^j$ in `$\tilde r$' gives $f_b^j=\sum_{i=0}^{Cj}\varphi_i\tilde r^i$, where $C=\kappa/\rho'$. Hence for the order of $\alpha_{tj}$ we find (ignoring the fact that we should reduce the coefficients of the product modulo $\tilde r$ at most once)
\[\ord_2(\alpha_{tj})\geq \delta(j+\min_{k=0}^{Cj}|t-k|).\]
It can then readily be checked that $\ord_2(\alpha_{tj})\geq\delta(|t|+(1-C)j)$ for $C<1$, and if $C\geq 1$ (in fact this could be forced) we distinguish between $t\geq 2Cj$, with $\ord_2(\alpha_{tj})\geq\delta(\frac 12 |t|+Cj)$ and $t\leq 2Cj$, where $\ord_2(\alpha_{tj})\geq \frac \delta{2C+1}(j+|t|)$. For proving that $g$, coming from the exact differential $dg$, can also be chosen in $T$, we need similar estimates using the full form of congruence (\ref{eq:congruence}). This congruence reads $(\ref{eq:congruence})=$
\[\frac 12d\left(\frac{X^j}3(4f+h^2)(2Y+h)\right)-d\int \left[\frac{X^j}2h(2f'+hh')+\frac j6X^{j-1}h(4f+h^2)\right]dX,\]
as can be checked by using the equality $(2Y+h)^2=4f+h^2$.\ep

\subsection{The differential equation.}\label{ssec:diffEq}

The goal of this section is to find the following commutative diagram:
\begin{equation}\label{eq:diagram}\begin{CD}
H_{MW}^- @>{\nabla}>> H_{MW}^-d\G\\
@VV{F_2}V @VV{F_2}V\\
H_{MW}^- @>{\nabla}>> H_{MW}^-d\G.
\end{CD}\end{equation}

Let us start with the definition of the connection $\nabla:H_{MW}\to H_{MW}d\G:t\mapsto \frac{\partial t}{\partial\G}d\G$. Similar computations as in the case of the differential $d$ show that $\frac{\partial}{\partial \G}$ and $\nabla$ are well defined on $T$ respectively $H_{MW}^\pm$. The expression for $\frac{\partial Y}{\partial\G}=\dot{Y}$ is similar to formula (\ref{eq:derivativeY}) where $'$ is replaced by $\dot{\ }$.

The map $F_2:T\to T$ represents a lift of the Frobenius automorphism $x\mapsto x^2$ in characteristic 2, and is defined\footnote{In fact $F_2$ equals $\sigma$ on $T$, but we prefer the notation $F_2$ to be used for the big modules.} as $\sigma$ on $\Qq$, $\G\mapsto\G^2$, $X\mapsto X^2$ and $Y$ maps to the unique solution in $T$ of $F_2(Y)^2+h^\sigma F_2(Y)-f^\sigma=0$ that is congruent to $Y^2$ modulo 2. It will follow from proposition \ref{prop:ConvergenceF} that with this definition $F_2(Y)$ actually sits in $T$. By extending $F_2$ with $dX\mapsto 2XdX$ and $d\G\mapsto 2\G d\G$ combined with the following lemma we have the two maps $F_2$ from the diagram above.
\begin{lemma}\label{lem:FpOnHMWminus}
The sum $\imath(F_2(YdX))+F_2(YdX)$ is exact.
\end{lemma}
\textsc{Proof.} Our proof is rather technical, we will use the sequence $W_k$ from the Newton iteration as in \cite{DenefVercauteren}, for which the approximation $F_2(Y)\equiv W_k\mod 2^k$ holds. Remark that this implies that $F_2(YdX)\equiv 2XW_kdX\mod 2^k$. We will show inductively for $k\geq 2$ that
\[\imath(W_kdX)+W_kdX\equiv 2^{k-2}h^2dX \mod 2^{k-1}.\]
As $W_2\equiv (f^\sigma-f^2-h^\sigma f)/h^2+(h^\sigma+2f)y/h\mod 2^2$ we find that $\imath(W_2dX)+W_2dX\equiv h^\sigma dX\mod 2$, and as $h^\sigma\equiv h^2\mod 2$ this satisfies our relation. Now the iterative step reads
\[h^2W_{k+1}\equiv -W_k^2+(h^2-h^\sigma)W_k+f^\sigma\mod 2^{k+1}.\]
Computing $h^2(\imath(W_{k+1})+W_{k+1})dX \mod 2^{k}$ and using the fact that $W_k^2\equiv f^\sigma-h^\sigma W_k\mod 2^k$ yields the equivalence
\[h^2(\imath(W_{k+1})+W_{k+1})dX \equiv 2h^2(\imath(W_kdX)+W_kdX)\mod 2^{k},\]
which in turn gives our induction.

It is possible to prove this lemma on a more conceptual level in the following sense: lifting from the coordinate ring of the curve in characteristic 2 to the Monsky-Washnitzer cohomology is functorial, and as Frobenius commutes with the involution below, it will also commute in the characteristic zero case.\ep

Finally we have that the above diagram is commutative, which can be seen for example by looking at the action of Frobenius and $\nabla$ on power series. We can derive from this the central differential equation. Let $F(\G)$ be the matrix of the operator $F_2$ on $H_{MW}^-$, given by $F_2(b_i)=\sum_kF_{ik}b_k$, and analogously let $G(\G)$ be the matrix of $\nabla$. Using the relation $\nabla\circ F_2=F_2\circ \nabla$ on basis elements the following equation is easily obtained:
\begin{equation}\label{eq:diffEqOrig}\dot F(\G)+F(\G)G(\G) = 2 \G G^{\sigma}(\G^2) F(\G).\end{equation}
We will come back later to the problem of solving this equation in a decent way.

Suppose that we use the same lift to some $\Q_{q^n}$ (including $\G\leftarrow\g$) in the algorithm of Denef and Vercauteren as we did here, then it is clear that if $F(0)$ equals their Frobenius in $\G=0$, the same will hold for $F(\g)$ for every $\g\in\S$ as $F(\G)$ is uniquely determined by (\ref{eq:diffEqOrig}) and $F(0)$.

\section{Behavior of matrices}\label{sec:matrices}
The theory in the foregoing section shows that the matrix of Frobenius $F(\g)$ for some $\g\in\S$ can be computed by working over a small field (for finding $F(0)$) and solving the right differential equation. One way to do this is by first finding $G$ and then using a recursive computation from equation (\ref{eq:diffEqOrig}). However --- as a general entry of $G$ is not a polynomial in $\G$ but rather a power series --- this would be rather slow, and in surplus we would need an expansion of $\nabla(b_i)$ which would require $\O(n^3)$ of memory. This section shows how to deal with this problem, and also gives an important estimate on $F(\G)$.

\subsection{Changing the matrices into some smaller form.}\label{ssec:ComputeSmallMatrices}

Define $v:=4f+h^2$ and $u := v'/2=2f'+hh'$. We construct a new basis for $H_{MW}^-$ as $d_i := vb_i$; the fact that this is a basis follows from proposition \ref{prop:DetBUnit}. The idea is that --- as $v$ arises as denominator in $\nabla b_i$ --- the basis $\{d_i\}$ gives in some sense a nicer matrix for the connection. We have (by definition) the following matrices, where the right hand sides are obtained by reduction using formula (\ref{eq:congruence}). By $(b_i)$ we mean a column vector of length $2g$ with $b_0$ on top.
\begin{align*}
\left(d_i\right)&\equiv B\cdot \left(b_i\right),\\
\nabla \left(b_i\right)&\equiv G\cdot \left(b_i\right),\\
\nabla \left(d_i\right)&\equiv D\cdot \left(b_i\right).
\end{align*}
As follows from the preceding section the entries of $G$ are elements of $S$, and it is not hard to see that the entries of $B$ and $D$ are polynomials in $\G$ over $\Qq$. Using these relations and the equality $\nabla\circ d=d\circ \nabla$ we find
\[D\cdot\left(b_i\right)\equiv\nabla\left(d_i\right)\equiv \dot{B}\cdot \left(b_i\right)+B\cdot\nabla\left(b_i\right)\equiv \dot{B}\cdot\left(b_i\right)+B\cdot G\cdot \left(b_i\right)\]
or in conclusion $D=\dot{B}+B\cdot G$.\\

\subsection{Adaptation of the differential equation}\label{ssec:newDiffEq}

If we combine the formula $D=\dot B+BG$ with the differential equation, we can find an equivalent equation where only polynomials of bounded degree --- see lemma \ref{lem:EstimatesBD} --- appear. We can however even go further, namely we will argue later on that we need in fact $r(\G)^MF(\G)$ for some positive integer $M$. Knowing $M$ we can find one `small' equation which has as solution precisely $K=r^MFB^{-1}$ and boundary or starting condition $K(0)=K_0$ for some relevant $K_0$.

We start with $\dot F(\G)+F(\G)G(\G) = 2 \G G^{\sigma}(\G^2) F(\G)$, hence multiplying with $B^{\sigma}$ on the left will remove $G^\sigma$. We suppress from now on the $\G$ and $\G^2$ from the equations.
\[B^{\sigma}\dot F+B^{\sigma}FG = 2 \G (D-\dot B)^{\sigma} F.\]
Next we substitute $KB=r^MF$, which after multiplying with $r^{M+1}$ leads to
\begin{equation}\label{eq:diffEqFinal}
(rB^{\sigma})\dot KB+(rB^{\sigma})KD+(-M\dot rB^\sigma+2\G r(\dot B-D)^\sigma)KB=0.\end{equation}
An important property of this equation is that all coefficients consist of polynomials of low degree. As proposition \ref{prop:DetBUnit} will show $B(0)$ is invertible, and hence it is possible to solve (\ref{eq:diffEqFinal}) using induction: write $K=\sum_i K_i\G^i$, where $K_0$ is known. Then we can find each $K_{k+1}$ one by one from $K_k, K_{k-1}, \ldots$ by looking at the coefficient of $\G^{k-1}$. Finally $r^MF$ is recovered as $KB$.

\subsection{Behavior of $B$ and $D$.}\label{ssec:BehaviorBD}

\begin{lemma}\label{lem:EstimatesBD}
For every $i,j$ there holds $\deg_{\G}B_{ij}\leq (2g+2)\kappa$ and $\ord_2(B_{ij})\geq -(3+\lfloor\log_2(5g+1)\rfloor)$ on the one hand, and on the other hand $\deg_{\G}D_{ij}\leq (2g+1)\kappa-1$ and $\ord_2(D_{ij})\geq -(3+\lfloor\log_2(5g)\rfloor)$.
\end{lemma}
\textsc{Proof.} We have for every $i$ the equivalence
\[(4f+h^2)X^iYdX\equiv \sum_{j=0}^{2g-1}B_{ij}X^jYdX.\]
The reduction formula (\ref{eq:congruence}) has to be applied at most $2g+1$ times, and each time $\deg_{\G}$ increases at most with $\kappa$. Bounding the denominator naively would give the following product which we give a name to be used in the next proposition,
\begin{equation*}\label{eq:productOrder}
P := \prod_{m=0}^{2g}\left(2(2g+1)+\frac{4m}3\right),\end{equation*}
which has order exactly $2g+1$. However, using lemma 2 of \cite{DenefVercauteren} gives the better logarithmic bound mentioned above. The results for $D$ can be proven with similar estimates.\ep

\begin{proposition}\label{prop:DetBUnit}
For every $\g\in\S$ we have $\ord_2(\det(B(\g)))=0$.
\end{proposition}
\textsc{Proof.} We will prove in a first step that
\[\det(B)\cdot P = \Res_X(u,v),\]
and afterwards some property of the resultant will show that for every $\g\in\S$ this last resultant has the same order $2g+1$ as $P$, which gives the proposition.\\

Define $\alpha_j := X^ju+(j/3)X^{j-1}v$ for $j := 0\ldots 2g$, then formula (\ref{eq:congruence}) reads $\alpha_jYdX\equiv 0$. We will suppress $YdX$ from the expressions during this proof, as they only make notation heavier. We define a square matrix $M$ over $\Q_q[\G]$ of dimension $4g+1$ which will be represented as a polynomial with coefficients in $\Qq[\G]$ and variables $\mu_0,\ldots,\mu_{2g},\lambda_0,\ldots,\lambda_{2g-1}$ and $X$. It has total degree 1 in the set of variables $\{\mu_i,\lambda_i\}$ and degree $4g$ in $X$. The entries of the matrix $M$ are given by the coefficients of $\mu_iX^j$ and $\lambda_iX^j$, enumerated in such a way that the first $2g+1$ rows correspond to $\mu_{2g}\ldots\mu_0$, the next rows to $\lambda_{2g-1}\ldots\lambda_0$, and the columns correspond to decreasing degrees of $X$. For example, the lower right entry is the coefficient of $\lambda_0X^0$. We start with
\[\lambda_0X^0v+\mu_0\alpha_0+\ldots+\lambda_{2g-1}X^{2g-1}v+\mu_{2g-1}\alpha_{2g-1} +\mu_{2g}\alpha_{2g}.\]
By means of the transformation $\lambda_j\leftarrow\lambda_j-((j+1)/3)\mu_{j+1}$ it is easy to see that the determinant of $M$ is precisely the resultant $\Res_X(u,v)$.

The reduction process gives rise to formulae of the form
\[X^jv=B_j(X)+\sum_{i=0}^{j+1}\beta_{ij}\alpha_i,\]
with $j=0\ldots 2g$, $\beta_{ij}\in\Qq[\G]$ and $\deg_X B_j(X)\leq 2g-1$. The coefficients of the $B_j$ are exactly the entries of the matrix $B$. If we substitute these expressions in our polynomial, we find
\[\lambda_0B_0+\ldots+\lambda_{2g-1}B_{2g-1}+ \left(\mu_0+\sum_{j=0}^{2g-1}\lambda_j\beta_{0j}\right)\alpha_0+\]\[ \left(\mu_{1}+\sum_{j=0}^{2g-1}\lambda_j\beta_{1j}\right)\alpha_{1}+\ldots+
\left(\mu_{2g}+\sum_{j=2g-1}^{2g-1}\lambda_j\beta_{2g,j}\right)\alpha_{2g}.\]
With the substitution $\mu_i\leftarrow\mu_i+\sum_{j=\max(i-1,0)}^{2g-1}\lambda_j\beta_{ij}$ again the determinant doesn't change, and the result is
\[\lambda_0B_0+\ldots+\lambda_{2g-1}B_{2g-1}+\mu_0\alpha_0+ \ldots+\mu_{2g}\alpha_{2g}.\]
In this form the upper half of the matrix is in `uppertriangular form', with $P$ as product of the elements on the diagonal. The lower half of the matrix has on the left only zeroes, and on the right the matrix $B$ appears (turned upside down and from left to right). This concludes the first part of the proof.

\begin{lemma}\label{lem:propertyResX}
Let $R$ be a ring and $\alpha,\beta,\gamma\in R[X]$ with $\deg\beta=\deg(\beta+\alpha\gamma)$, then $\Res_X(\alpha,\beta)=\Res_X(\alpha,\beta+\alpha\gamma)$.
\end{lemma}
This lemma remains true without the condition on the degree, given that $\alpha$ is monic. Otherwise the resultants agree up to an appropriate power of the leading coefficient of $\alpha$.\\
\textsc{Proof.} The matrix defining the second resultant can be achieved from the matrix defining the first resultant by adding to the rows according to $\beta$ suitable multiples of the rows of $\alpha$. These elementary row operations do not change the determinant.\ep

\noindent We write $\Res_X(v,u)=\Res_X(H,2f'+hh')\cdot\Res_X(4Q_f+(h^2/H),2f'+hh')$. By the lemma and the fact that $H$ and $Q_fH'$ are relatively prime we have that the first factor has order $\deg H$. Define $\tilde{h} := h/H$, then we have --- as can be checked by writing $\tilde{h}$ in a product of linear factors over $\Q_q^{\text{alg cl}}$ of $H$ --- that $\tilde{h}$ is a divisor of $H\tilde{h}'$ with integral quotient $\alpha$. The lemma implies that
\[\Res_X(4Q_f+\tilde{h}h,2f'+hh')=\Res_X(4Q_f+\tilde{h}h,2f'+hh'- (H'+\alpha)(4Q_f+\tilde{h}h))\]\[ = \Res_X(4Q_f+\tilde{h}h,2Q_f'H-2Q_fH'-4Q_f\alpha).\]
Remark that the coefficient of $X^{2g}$ of the second polynomial in these equalities is always congruent to 2 modulo 4, and hence nonzero.

The last resultant above equals $2^{\deg Q_f}$ times
\[\Res_X(4Q_f+\tilde{h}h,Q_f'H-Q_fH'-2Q_f\alpha),\]
and reducing this result modulo 2 gives $\Res_X(\tilde{h}h,Q_f'H-Q_fH')$. Again using the lemma we find $\Res_X(\tilde{h}h,-Q_fH')$ modulo 2, which is nonzero by construction. In conclusion we see that $\Res_X(v,u)$ has an order of exactly $\deg Q_f+\deg H=2g+1$.\ep

A consequence of this proposition is an estimate on $B^{-1}$. Indeed, suppose $2^\varepsilon B$ is integral, then the fact that the inverse of a matrix equals its adjunct matrix divided by the determinant gives that the order of $B^{-1}$ is at least $-(2g-1)\varepsilon$. Together with lemma \ref{lem:EstimatesBD} we can conclude that, defining $\beta':=(2g-1)(3+\lfloor\log_2(5g+1)\rfloor)=\O(g\log g)$, we have $\ord_2(B^{-1})\geq -\beta'$.

\subsection{On the convergence rate of $F(\G)$.}\label{ssec:ConvergenceF}

\begin{proposition}\label{prop:ConvergenceF}
Let $N\in\N$ and $f(\G)$ be an entry of $F(\G)$, reduced modulo $2^N$. Then there exist explicit constants $\chi_1=\O(N\tilde D)$ and $\chi_2=\O(g\kappa N\tilde D)$ such that $r^{\chi_1}f(\G)$ is a polynomial of degree at most $\chi_2$. Also we have an explicit constant $\varphi=\O(\log g)$ such that $\ord_2F(\G)\geq -\varphi$.
\end{proposition}
\textsc{Proof.} Recall from \cite{DenefVercauteren} the approximation $W_k$ to $F_2(Y)$, also used in the proof of lemma \ref{lem:FpOnHMWminus}. By defining $\alpha_k(X,\G)$, $\beta_k(X,\G)$ such that $W_k=\alpha_k+Y\beta_k$; $\Delta_{\alpha,k}:=(\alpha_k-\alpha_{k-1})/2^{k-1}$ and similar $\Delta_{\beta,k}$ we can compute $H^{2\tilde D}W_{k}$ from the following formula of \cite{DenefVercauteren}.
\[H^{2\tilde D}W_k\equiv Q_h^2\cdot\left\{ -W_{k-1}^2+(h^2-h^\sigma)W_{k-1}+f^\sigma \right\}\mod 2^k.\]
This gives as result:
\[H^{2\tilde D}W_k\equiv -Q_h^2\sum_{1\leq i<j, i+j\leq k}2^{i+j-1}\left( \Delta_{\alpha,i}\Delta_{\alpha,j}+(f-hY)\Delta_{\beta,i}\Delta_{\beta,j} \right)\]
\[-YQ_h^2\sum_{i+j\leq k}2^{i+j-1}\Delta_{\alpha,i}\Delta_{\beta,j}- Q_h^2\sum_{2i\leq k+1}2^{2(i-1)}\left(\Delta_{\alpha,i}^2+(f-hY)\Delta_{\beta,i}^2\right)\]
\[+(h^2-h^\sigma)Q_h^2\sum_{i\leq k-1}2^{i-1} \left(\Delta_{\alpha,i}+\Delta_{\beta,i}Y\right)+Q_h^2f^\sigma\mod 2^k.\]
We start by proving that the numerators in $W_k$ --- the right hand side of the equation above, expanded as an $H$-adic series --- for $k\geq 2$ have $\deg_\G$ at most $Ak-B$, with $\delta:=\omega-\kappa$, $A:=\omega+\delta$, $B:=A+\delta$, and
\[\omega := 2\kappa+\deg_\G Q_h^2+[(\deg_Xf^2+2\deg_XQ_h)/\deg_XH+3]\eta.\]
Here $\omega=2A-B$ is a bound for the degree in $\G$ in $W_2$, as can be checked by an easy computation. To prove the bound $Ak-B$ we use induction and consider each term in the formula for $W_k$ above, for instance for $Q_h^2\Delta_{\alpha,i}\Delta_{\alpha,j}$ with $i\geq 2$, $j\geq 2$ we find as bound
\[Ai-B+Aj-B+(\deg_XQ_h^2/\deg_XH+2)\eta\leq Ak-2B+(\ldots)\eta\leq Ak-B.\]
As $A$ is also a bound for the numerators in $W_1$ and $i,j\leq k-1$ we have our estimate for all $i,j$. The term with $\eta$ comes from expanding polynomials in $X$ as series in $H$.\footnote{The case where $H$ is a constant is similar but easier.} For the other terms a similar computation works, for example for $Q_h^2f\Delta_{\beta,i}^2$ we have, as $2i\leq k+1$,
\[2Ai-2B+\kappa+[(\deg_Xf+2\deg_XQ_h)/\deg_XH+3]\eta\leq Ak-B.\]

In a second step we have to reduce $W_k$ in the cohomology. As $F_2(Y)\in H_{MW}^-$ we can confine ourselves to the part with $Y$ in it. First we take some $g(\G)X^rYdX$, and reducing this by using formula (\ref{eq:congruence}) adds less than $r\kappa$ as degree in $\G$. Lemma 1 of \cite{DenefVercauteren} shows that $X^r$ has possible nonzero coefficient modulo $2^M$ only if $r\leq (aM+b)s$ with $as=2(2g+1-2\deg_Xh)$ and $bs=7\deg_Xh-3(2g+1)$. Take $M$ such that $M-(3+\log_2((aM+b)s+g+1))\geq N$, then clearly $M=\O(N)$ and lemma 2 of \cite{DenefVercauteren} gives that it is enough to compute $W_M$ for finding $F_2(Y)\mod 2^N$ in $H_{MW}^-$, at least for the part without denominators $H$. Thus the worst possible $\deg_{\G}$ comes from the term $VH^{aM+b}$, which has $\deg_{\G}$ as most $AM-B$. During the reduction an extra $(aM+b)s\kappa$ can occur, and in conclusion the contribution of the part without denominator $H$ is at most $AM-B+(aM+b)s\kappa$.

For the second part of $F_2(Y)$ we consider terms of the form $V/H^\ell YdX$ for $\ell>0$. During the reduction from $1/H^\ell$ to $1/H^{\ell-1}$ the degree in $X$ increases with at most $s+2g$ and the degree in $\G$ with at most $2g\kappa$. Also a denominator $r(\G)$ appears. In the end we also have to reduce as in the previous paragraph, starting from $\deg_X$ at most $\ell(s+2g)$. Let $\tilde{a}:=4\tilde D$ and $\tilde{b}:=-6\tilde D$, so that lemma 1 of \cite{DenefVercauteren} implies that modulo $2^{\tilde{M}}$ we only need $\ell\leq \tilde{a}\tilde{M}+\tilde{b}$. Then with $\tilde{M}$ such that $\tilde{M}-(3+\log_2(\tilde{M}+1))\geq N$ again $\tilde{M}=\O(N)$ and from lemma 3 of \cite{DenefVercauteren} it follows that $W_{\tilde{M}}$ suffices for this part. Hence the worst case here is the denominator $H^{\tilde{a}\tilde{M}+\tilde{b}}$, where $\deg_\G$ is at most $A\tilde{M}-B$. All together this gives a degree in $\G$ of at most
$A\tilde{M}-B+2g\kappa(\tilde{a}\tilde{M}+\tilde{b})+ (\tilde{a}\tilde{M}+\tilde{b})(s+2g)\kappa$, and a denominator $r^{\tilde{a}\tilde{M}+\tilde{b}}$.

It is now easy to find the bounds from the lemma: the denominator is $r^{\tilde{a}\tilde{M}+\tilde{b}}$ with $\tilde{a}=\tilde{b}=\O(\tilde D)$ and $\tilde{M}=\O(N)$; and as bound for the degree we find
\[\max\left\{A\tilde{M}-B+2g\kappa(\tilde{a}\tilde{M}+\tilde{b})+ (\tilde{a}\tilde{M}+\tilde{b})(s+2g)\kappa,AM-B+(aM+b)s\kappa\right\}.\]
Using $A=B=\O(g\kappa)$, $s=\O(g)$ and $as$ and $bs$ as before the lemma follows.

Remark that we should in fact look at $F_2(X^iY)$ for $i=0\ldots 2g-1$, but the possible increased $\deg_{\G}$ caused by this is absorbed in the rough estimates during the proof.\\

In order to determine $\varphi$ we need to combine lemmata 1, 2 and 3 of \cite{DenefVercauteren}. Choosing a modulus $2^k$ lemma 1 implies that the highest appearing degree of $X$ in the $Y$-part of $F_2(Y)$ is less than $(4g+2)k+g$. Linked with lemma 2 this part gives then an order bigger than
\begin{equation}\label{eq:bound1FY}
\min_{k\geq 0}\left(k-3-\log_2((4g+2)k+2g+1)\right).
\end{equation}
On the side with denominators we find as extremum $4\tilde Dk-6\tilde D$, and lemma 3 then gives the bound
\begin{equation}\label{eq:bound2FY}
\min_{k\geq 0}\left(k-3-\log_2(4\tilde Dk-6\tilde D+1)\right).
\end{equation}
Now we can take $-\varphi$ as the minimum of (\ref{eq:bound1FY}) and (\ref{eq:bound2FY}), and we see immediately that $\varphi=\O(\log g)$.\ep

\begin{note}\label{note:orders} When implementing these results one finds that $F(\G)^{-1}$, the matrix of the big Frobenius and $B^{-1}$ actually have also very good $2$-adic valuation\footnote{This is also true for $F(\G)^{-1}$ and the big Frobenius in odd characteristic.}, good enough to suggest a bound of $\O(\log g)$ for them as well. However, we do not know of a way to prove this, but in \cite{DenefVercauteren} a heuristic argument is given for the big Frobenius. It is worth noting that a proof of these results would diminish all our complexity estimates by a factor $g$.\end{note}

\subsection{Error propagation in the inductive computation.}\label{ssec:errorProp}

When solving the equation \begin{equation}\label{eq:diffEqFinal1}
(rB^{\sigma})\dot KB+(rB^{\sigma})KD+(-M\dot rB^\sigma+2\G r(\dot B-D)^\sigma)KB=0,\ \ K(0)=K_0\end{equation}
in an inductive manner, we could estimate the loss in accuracy in a naive way. However, already $\dot K=\sum iK_i\G^{i-1}$ implies division by $k$ for computing $K_k$, and hence at least $\ord_2((N_{\G}-1)!)$ would be lost as accuracy, assuming working modulo $\G^{N_{\G}}$. It turns out to be possible to do better, as we will show in theorem \ref{thm:errorProp}. Some form of this theorem has been found independently from the author by Gerkmann \cite{GerkmannEC}.

Let $-\varphi$ be the lower bound for the order of $F(\G)$ found above, and $-\varphi_0$ a bound for $F(\G)^{-1}$. By lemma 24 in \cite{HubrechtsHECOdd} --- the proof of which is also correct for $p=2$ --- we can take $\varphi_0= \varphi(2g-1)+g$. Denote with $\mathcal{K}$ the solution of (\ref{eq:diffEqFinal1}) obtained by working modulo $2^N$ and starting with $\mathcal{K}_0=K_0=r(0)^MF_0B_0^{-1}$. $K$ itself will denote the exact solution, hence $K=r^MFB^{-1}$. Finally we write $A_0$ for $r(0)^MF_0=K_0B_0$.

\begin{theorem}\label{thm:errorProp}
With $\tilde K:=2^{-N}(\mathcal{K}-K)= \sum_i \tilde{K}_i\G^i$ we have
\[\ord_2(\tilde{K}_i)\geq -(2g\varphi + g+1)\cdot \log_2(i+1)-\alpha,\]
where $\alpha:=(12g-1)(3+\lfloor\log_2( 5g+1)\rfloor)+(10g-1)\varphi+5g$.
\end{theorem}
\textsc{Proof.} We will prove this theorem in a number of steps. Let us first define and recall some terms. For ease of notation we write $E := -M\dot rB^\sigma+2\G r(\dot B-D)^\sigma$. We know the following bounds:
\[\ord_2(B)= \ord_2(B^\sigma)\geq -\beta:=-(3+\lfloor\log_2( 5g+1)\rfloor),\]
and with $\beta'$ such that $\ord_2(B^{-1})=\ord_2((B^\sigma)^{-1})\geq -\beta'$ as defined after the proof of proposition \ref{prop:DetBUnit} we have $\beta+\beta'=2g\beta$. The same way we have $\varphi+\varphi_0=2g\varphi +g$.
\begin{definition}\label{def:LogConvergence}
Let $A_i$ be for every $i\geq 0$ a $(d\times d)-matrix$ over $\C_2$ and $x,y\in\mathbb{R}$. We say that a power series $\sum_i{A_i\G^i}$ converges $(x,y)$-logarithmically if for all i
\[\ord_2{A_i}\geq -x\log_2(i+1) - y.\]\end{definition}
To shorten notation we will often forget the word  `logarithmically'.
\begin{lemma}\label{lem:prodLogConv}
If $\sum_iA_i\G^i$ and $\sum_iB_i\G^i$ converge $(x,y)$- respectively $(x',y')$- logarithmically, then their product converges as $(\max(x,x'),y+y')$.\end{lemma}
\textsc{Proof.} The coefficient of $\G^k$ in the product is $\sum A_iB_j$, summed over $i+j=k$. Hence the $2$-order is at least $$-x\log_2(i+1)-x'\log_2(j+1)-(y+y'),$$
and as $\log_2(k+1)\leq\log_2(i+1)+\log_2(j+1),$ we find the lemma.\ep
\vspace{-\baselineskip}
\begin{lemma}\label{lem:ConvergenceOfC}
Let $C$ be the (exact) solution of $\dot C B+CD=0$ subject to $C(0)=B_0^{-1}$, then $C$ converges $(\varphi+\varphi_0,\beta')$-logarithmically, and for $C^{-1}$ with $C^{-1}(0)=B_0$ we find $(\varphi+\varphi_0,\beta)$.\end{lemma}
\textsc{Proof.} The matrix $C' := CB$ gives in fact the solutions around zero of the equation $\nabla=0$, or $\dot C'+C'G=0$, and from the diagram (\ref{eq:diagram}) we can deduce the equality
\[C'^{\sigma}(\G^2)F(\G)=F(0)C'(\G).\]
Now exactly the same proof as for proposition 25 in \cite{HubrechtsHECOdd} gives that $C'$ converges as $(\varphi+\varphi_0,0)$. As $B^{-1}$ can be considered to converge as $(0,\beta')$, lemma \ref{lem:prodLogConv} gives the result. The estimate for $C^{-1}=B(C')^{-1}$ can be proven in a similar fashion.\ep

We now give an estimate on the error propagation for two `partial solutions' of the equation. Remark that we don't need these in the algorithm, only in this proof. A lemma with the flavor of the following one was first given by Lauder \cite{LauderRecursive}, but we give a proof similar to our proof in \cite{HubrechtsHECOdd}. Let $\mathcal{C}$ be the solution computed inductively modulo $2^N$ from the equation $\dot CB+CD=0$ with $\mathcal{C}(0)=B_0^{-1}$.
\begin{lemma}\label{lem:errorPropC}
$2^{-N}(\mathcal{C}-C)$ converges $(\varphi+\varphi_0+1,\beta+2\beta')$-logarithmically.\end{lemma}
\textsc{Proof.} It is easy to see (a formal argument will be given later) that $\mathcal{C}$ satisfies $\dot{\mathcal{C}}B+\mathcal{C}D=2^N\mathcal{E}_1$ with $\mathcal{E}_1$ some matrix of power series with integral coefficients. Let $L$ be such that $2^NLC=\mathcal{C}-C$. Then we have the equalities $$2^N\mathcal{E}_1=\dot{\mathcal{C}}B+\mathcal{C}D-\dot CB-CD =2^N(\dot LCB+L\dot CB+LCD)=2^N\dot LCB$$ and as a consequence $\dot L=\mathcal{E}_1B^{-1}C^{-1}$. If we integrate $\dot L$ we find as integration constant $L_0=0$, and hence
\[2^{-N}(\mathcal{C}-C) = LC = \left(\int \mathcal{E}_1B^{-1}C^{-1}d\G\right) C.\]
As integrating is not worse then adding 1 to the logarithmic factor, we find the lemma.\ep

Let $\mathcal{P}$ and $P$ be the computed modulo $2^N$ resp.\ exact solution of $(rB^\sigma)\dot P+EP=0$ subject to $P(0)=I$, then a trivial computation shows that $K=PA_0C$ satisfies (\ref{eq:diffEqFinal1}). Now the previous results give that $P=KC^{-1}A_0^{-1}$ converges $(\varphi+\varphi_0,\beta+\beta'+\varphi+\varphi_0)$-logarithmically  and the same holds for $P^{-1}=A_0CK^{-1}$. Working similar to lemma \ref{lem:errorPropC} shows that $2^{-N}(\mathcal{P}-P)$ converges as $(\varphi+\varphi_0+1,2\beta+3\beta'+2(\varphi+\varphi_0))$. In this we use $(rB^\sigma)\dot{\mathcal{P}}+E\mathcal{P}=2^N\mathcal{E}_2$.

The proof of the theorem can now be completed by estimating $\mathcal{K}-\mathcal{P}A_0\mathcal{C}$ and $\mathcal{P}A_0\mathcal{C}-K$ and summing these terms. Denote the additive operator of (\ref{eq:diffEqFinal1}) by $\Delta$, hence (\ref{eq:diffEqFinal1}) equals $\Delta K=0$.
\begin{lemma}\label{lem:errorPropFirstPart}
$2^{-N}(\mathcal{K}-\mathcal{P}A_0\mathcal{C})$ converges $(\varphi+\varphi_0+1,5\beta+6\beta'+5\varphi+4\varphi_0)$-logarithmically.
\end{lemma}
We will first show how to see that solving $\Delta K=0$ inductively modulo $2^N$ amounts to $\Delta\mathcal{K}=2^N\mathcal{E}$ for some integral matrix $\mathcal{E}$. For each $k$ we compute $\mathcal{K}_k$ from
\[\left[r(0)B^\sigma_0\mathcal{K}_kB_0+f_k(\mathcal{K}_{k-1},\mathcal{K}_{k-2},\ldots) \right]\G^{k-1}=2^N(\text{integral error matrix})\G^{k-1}\]
for some linear functions $f_k$. The sum over all these equations gives $\Delta\mathcal{K}=2^N\mathcal{E}$.

Let $L$ be defined such that $2^NPLA_0C=\mathcal{K}-\mathcal{P}A_0\mathcal{C}$, then we compute
\begin{equation}\label{eq:DeltaEqn} 2^{-N}(\Delta\mathcal{K}-\Delta(\mathcal{P}A_0\mathcal{C}))=\Delta (PLA_0C)=rB^{\sigma}P\dot LA_0CB.\end{equation}
Using the same integral as before and the fact that
\[\Delta(\mathcal{P}A_0\mathcal{C})= 2^N(rB^\sigma\mathcal{P}A_0\mathcal{E}_1+ \mathcal{E}_2A_0\mathcal{C}B),\]
we find our result. Indeed, for $2^{-N}\Delta(\mathcal{P}A_0\mathcal{C})$ we find $(\varphi+\varphi_0, 2\beta+\beta'+2\varphi+\varphi_0)$, and adding the inverse of the factors in the right hand side of (\ref{eq:DeltaEqn}) gives the lemma.\ep

To control the difference $2^{-N}(\mathcal{P}A_0\mathcal{C}-PA_0C)$ we add a cross term:
\[2^{-N}(\mathcal{P}A_0\mathcal{C}-PA_0\mathcal{C}+PA_0\mathcal{C}-PA_0C)= 2^{-N}(\mathcal{P}-P)A_0\mathcal C+2^{-N}PA_0(\mathcal{C}-C).\]
The $(\varphi+\varphi_0+1,4\beta+5\beta'+3\varphi+2\varphi_0)$-logarithmic convergence of this difference is now clear, and taking the maximum of this result and the last lemma gives the theorem.\ep

\section{The algorithm}\label{sec:algorithm}
In this section we give a concrete presentation of the algorithm. We suppose that the polynomials $\bar{H}(X,\G)$, $\bar{h}(X,\G)$ and $\bar{f}(X,\G)$ are given as explained in section 2. The input for the algorithm is hence formed by these polynomials over $\Fq=\F_{2^a}$ and some allowable parameter $\og\in\F_{q^n}$. The output is the zeta function of complete model of the hyperelliptic curve given by $Y^2+\bar{h}(X,\og)Y=\bar{f}(X,\og)$.\\

\noindent \textsc{Step 1.} Compute $\Q_q$ as explained in section \ref{ssec:pAdicArith}, lift $\bar{H}$, $\bar{Q}_{\bar{f}}$ and hence also $\bar{h}$ and $\bar f$ to $\Qq$ such that $H$ and $Q_f$ remain monic, and compute the resultant $r(\G)=\Res_X(H,Q_f\cdot H')$. Let $g$ be the genus and $M=\chi_1$, $\chi_2$ and $\varphi$ as follows from the proof of proposition \ref{prop:ConvergenceF} with $N$ defined as below. Also $\varphi_0 := \varphi(2g-1)+g$. Define
\[N_f := \left\lceil\log_2{2g\choose g}+1+ang/2\right\rceil, \quad N := N_f + an\varphi + 2gan\varphi, \quad N_{\G} := \chi_2+1,\]
\[\qquad N_2 := N + 12g(3+\lfloor\log_2( 5g+1)\rfloor)+(10g-1)\varphi+5g.\]
From now on we work modulo $2^{N_2}$ (in the beginning of the algorithm) and $\G^{N_{\G}}$.\\
\textsc{Step 2.} Compute the matrices $B$ and $D$ by using formula (\ref{eq:congruence}).\\
\textsc{Step 3.} Calculate $F(0)$ as explained in \cite{DenefVercauteren}, but with the higher accuracy $2^{N_2}$. Remark that we need the \emph{small} Frobenius.\\
\textsc{Step 4.} Compute $K$ inductively from the equation (\ref{eq:diffEqFinal1}) with starting condition $K_0=r(0)^MF(0)B(0)^{-1}$, and find $F'(\G):=r(\G)^MF(\G)$.\\
After step 4 we can switch to the accuracy $N$ instead of $N_2$.\\
\textsc{Step 5.} Let $\bar{\psi}(z)$ be the minimal polynomial of $\og$ over $\F_q$, where we suppose\footnote{This is not crucial, if $\og$ defines a smaller field then the zeta function over $\F_{q^n}$ is easily derived from it.} that $\F_q(\og)=\F_{q^n}$, and let $\psi(z)$ be the Teichm\"{u}ller modulus lift of $\bar{\psi}$ as explained in section \ref{ssec:pAdicArith}. Then $\Q_{q^n}=\Q_q[z]/\psi(z)$ and $z$ is the Teichm\"{u}ller lift of $\og$. Determine $$F(z)=\frac 1{r(z)^{\chi_1}}\cdot F'(z).$$
\textsc{Step 6.} Compute $\mathcal{F}=\prod_{i=1}^nF(z)^{\sigma^{n-i}}$ as explained by Kedlaya in \cite{KedlayaCountingPoints} and find $Z(T)$ as the polynomial $\det(I-\mathcal FT)$ with coefficients between $-2^{N_f-1}$ and $2^{N_f-1}$. Output now $Z(T)\cdot[(1-T)(1-2^{an}T)]^{-1}$.

\begin{proposition}\label{prop:AlgorithmCorrect}
The above algorithm returns the correct result.
\end{proposition}
\textsc{Proof.} The Lefschetz fixed point formula on the Monsky-Washnitzer cohomology gives as explained in \cite{DenefVercauteren} that the result is correct if we can compute $\mathcal{F}$ and $Z(T)$ exactly, and the theory from section \ref{sec:analysis} and \ref{sec:matrices} implies that if every step was done with exact precision, we would indeed find the required matrix $\mathcal{F}$. As we cannot work with this infinite precision, we need to show that the chosen accuracy is high enough. From the Weil conjectures it follows that $\mathcal{F}\bmod 2^{N_f}$ is sufficient to recover the zeta function, and proposition \ref{prop:ConvergenceF} proves that $N_{\G}$ suffices to compute $r(\g)^MF(\g)$ modulo $2^N$. The crucial difficulty is to control the loss of precision introduced by working with non integral elements of $\Qq$. It is clear that computing $r$, $B$ and $D$ gives no significant loss in precision. For computing $KB$ we can bound the introduced error as in theorem \ref{thm:errorProp}. This gives that the loss in precision is at most $12g(3+\lfloor\log_2( 5g+1)\rfloor)+(10g-1)\varphi+5g$. Here we have added $\beta$ to the result of theorem \ref{thm:errorProp}.

We should also take notice of possible loss in accuracy in the computation of $\mathcal{F}$ as a product, which requires an extra $an\varphi$ of accuracy. But as pointed out in note \ref{note:orders}, in practice $\mathcal{F}$ turns out to have a similar order as $F(\g)$, hence this increment of $N$ can in practice be chosen lower. Another problem appears in the computation of the characteristic polynomial of $\mathcal{F}$. One naive way of doing this would be to compute the trace of $\mathcal{F}^i$ for $i=1\ldots 2g$ and to use Newton's formula
\[\det(I-\mathcal{F}T)=\exp\left(-\sum_{k=1}^\infty\text{Tr}(\mathcal{F}^k) \frac{T^k}k\right),\] which would require an extra precision of $2g+\log_2(2g)$ from the exponential and the denominators $k$, and $2gan\varphi$ for the trace of $\mathcal{F}^{2g}$. A better way however is explained in \cite{CastryckDenefVercauteren}. Here we first make $\mathcal{F}$ integral by multiplying it with some power of $2$, and then use a slightly altered version of reduction to the Hessenberg form of a matrix, suitable for working in $\Z_{q^n}$. The loss in precision is then $2gan\varphi$. We can conclude that the values of $N$ and $N_2$ are sufficient.

\section{Complexity analysis}\label{sec:complexity}
\subsection{$2$-Adic arithmetic}\label{ssec:pAdicArith}
As central source for this section we use chapter 12 by Vercauteren of \cite{HandboekHECC}, and we always assume asymptotically fast arithmetic, meaning that basic operations can all be done in essentially linear time. We suppose here that we are working modulo $2^N$, hence representing an element of $\Q_2$ takes $\O(N)$ bits (if its order is not too low) and computing with it $\Ot(N)$ bit operations. Remember that $q=2^a$. Let $\Fq\cong \F_2[x]/\bar{\chi}(x)$, then we define $\Qq\cong \Q_2[x]/\chi(x)$ where $\chi$ is the Teichm\"{u}ller modulus that projects to $\bar{\chi}$. A Teichm\"{u}ller modulus is a minimal polynomial for Teichm\"{u}ller lifts, or equivalently $\chi(x)|x^q-x$. In \cite{HandboekHECC} an algorithm of Harley is given that computes $\chi$ in time $\Ot(aN)$. Basic operations, including the $2$nd power Frobenius automorphism $\sigma$, need the same amount of time.

If $\bar{\psi}(z)$ is the minimal polynomial of $\bar{\gamma}$ over $\Fq$, we can compute the Teich\-m\"{u}l\-ler modulus $\psi(z)$ over $\Qq$ as follows. First determine $\varphi(y)$ such that $\Q_{q^n}\cong \Q_2[y]/\varphi(y)$, $\varphi(y)|y^{2^{an}}-y$ and $\bar\varphi(\og)=0$ as above, in time $\Ot(anN)$. Second, as $\varphi(z)=0$, we have that $\psi|\varphi$, or $\varphi=\psi\cdot \psi'$. Now $\bar{\psi}$ and $\bar{\varphi}$ are known, hence $\bar{\psi}'$ can be recovered easily, and using Hensel lifting as in \cite{ModernCompAlg} gives $\psi$ in time $\Ot(anN)$. Again this is also the time required for basic operations.

Computing $\sigma^k$ of an element of $\Q_{q^n}$ can be done trivially by applying $k$ times $\sigma$, resulting in a complexity of $\Ot(kanN)$. However, further on it will be advantageous to be able to compute $\sigma^k(z)$ in a faster way. Indeed, we can compute $\og^{2^k}$ in time $\Ot(kan)$ by repeated squaring, and using the generalized Newton lifting of \cite{HandboekHECC} we find the Teichm\"{u}ller lift of $\og^{2^k}$, which equals $\sigma^k(z)$, in time $\Ot(anN)$.

\subsection{Analysis of the algorithm}\label{ssec:alysisalgorithm}
We use the $2$-adic arithmetic always as in the previous paragraph. Let $\omega$ be an exponent for matrix multiplication, meaning that multiplying two $k\times k$ matrices over some ring $R$ takes $k^\omega$ operations in $R$. We can take $\omega=2,376$. It is easy to check the following bounds:
\begin{align*}
\varphi&=\O(\log g)=\Ot(1),\qquad \varphi_0=\O(g\log g)=\Ot(g),\\
N_f&=N=N_2=\O(ang\log g)=\Ot(ang),\\
N_{\Gamma}&=\Ot(g\k N\tilde D)=\Ot(g^2 a\k n\tilde D).
\end{align*}
Computing the lifts of $\bar{H}$ and $\bar{Q}_{\bar{f}}$ costs essentially nothing, and the computation of the resultant $r(\G)$ can be achieved in time $\Ot(g^{1+\omega} aNg\k)=\Ot(g^{3+\omega}a^2\k n)$, see e.g.\ \cite{Villard}, where we use the fact that we are working with polynomials in $\G$ of degree at most $\O(g\k)$. To determine $B$ and $D$ we have to use formula (\ref{eq:congruence}) at most $\O(g)$ times, and each step requires time $\Ot(aN\cdot g\k\cdot g)$, which comes from `$\Q_q\cdot \deg_{\G}\cdot\deg_X$'. Together this gives $\Ot(g^4a^2n\k)$. Next we have the recursive formula for finding $K$. Each of the $N_{\G}$ steps consists of $\O(g\k)$ multiplications of matrices whose entries have size $\O(aN)$, resulting in $\O(g\k g^\omega a N N_{\G})=\Ot(g^{4+\omega}a^3\k^2n^2\tilde D)$. The size of $K$ is $\O(g^2 aN N_{\G})=\Ot(g^5 a^3\k n^2\tilde D)$, which will be the overall memory requirements of the algorithm. Remark that we can ignore the operations for finding $B^\sigma$ and the like.

Repeating the complexity analysis of \cite{DenefVercauteren}\footnote{In that paper the memory requirements are actually $\log g$ bigger than written there, because the computation of the characteristic polynomial of the big Frobenius needs to take care of the emerging denominators. However, as we are only interested in the small Frobenius, this factor does not appear.}, we can confine ourselves to the worst case mentioned there, and as we skip the computation of the norm of the matrix, the most time consuming step is step 4 of the algorithm, which takes $\Ot(g^3aN^2)=\Ot(g^5a^3n^2)$. The memory requirements are $\O(g^4a^3n)$. The minimal polynomial $\bar{\psi}$ can be computed in time $\Ot(an\sqrt{an}+(an)^2)$, see \cite{ShoupMinimalPolynomial}, and finding $\psi$ out of $\bar{\psi}$ takes $\Ot(anN)$ bit operations.

Let $f(\G)$ be an entry of $r(\G)^{\chi_1}F(\G)$, then we need to find $f(z)$, a substitution $\G\leftarrow z$ that can be done very fast using our Teichm\"{u}ller modulus. Indeed, we just have to reduce $f(z)$ modulo $\psi(z)$, which takes for the whole of the matrix $\Ot(g^2anN_{\G})=\Ot(g^5a^2\k n^2\tilde D)$ bit operations. Division by $r(z)^{\chi_1}$ is again neglectable. Remark that until now, where we have found the matrix of the small Frobenius, our algorithm has complexity $\Ot(n^2)$ in $n$.

For the last step Kedlaya's method consists of the following iteration: $M_0:=F(z)$ and $M_{i+1}=M_i^{\sigma^{2^i}} M_i$. This requires $\log n$ times a matrix multiplication over $\Q_{q^n}$, which needs time $\Ot(g^{\omega} anN)$. The computation of $\sigma^k$ on $4g^2$ elements requires $\Ot(g^2\cdot(k=an)\cdot anN)$ bit operations.

Combining all these facts gives up to step 5 a complexity of $\Ot(g^{4+\omega}a^3\k^2n^2\tilde D)$ bit operations and $\Ot(g^5a^3\k n^2 \tilde D)$ bits of memory. Now as `on average' $\tilde D=\O(1)$ --- worst case being $\tilde D=\O(g)$ --- this gives the first term in the first complexity and the memory requirements in theorem \ref{thm:princThm}. Step 6 gives the second part of the time estimate.

\section{Improvements}\label{sec:conclusion}
\subsection{Subcubic counting.}\label{ssec:subcubic}
The most time consuming step in the above algorithm is in fact the determination of $F(z)^{\sigma^k}$ for $k$ of the order $\O(an)$, taking time $\Ot(ga^3n^3)$. It is however possible to do this with a faster method. Let $\alpha(z)\in\frac{\Q_q[z]}{\psi(z)}$, then the equality $\alpha(z)^{\sigma^k}=\alpha^{\sigma^{k \bmod a}}(z^{\sigma^k})$ shows that we only have to compute $4g^2\log n$ times $\alpha^{\sigma^\ell}(z^{\sigma^k})$ with $\ell=\O(a)$ and $k=\O(an)$, where $\alpha$ is a polynomial modulo $2^N$ over $\Qq$ of degree at most $n-1$. The computation of $\alpha^{\sigma^\ell}$ takes at most time $\Ot(aN\ell n)=\Ot(ga^2n^2)$. On the other hand we have the \emph{modular composition of polynomials} $\alpha^{\sigma^\ell}(z^{\sigma^k})$. As said before the computation of $z^{\sigma^k}$ takes only $\Ot(ga^2n^2)$ time, and as explained in \cite{HubrechtsHECOdd} this composition can be achieved in time $\Ot(ga^2n^{2,667})$, at the cost of an increase in memory use, resulting in $\O(ga^2n^{2,5})$. This proves theorem \ref{thm:subcubic} from the introduction.

\subsection{Lots of curves.}\label{ssec:lotscurves}
Using fast multipoint evaluation \cite{ModernCompAlg} it is possible to compute $\O(n)$ zeta functions within one family in time and memory usage $\Ot(n^3)$. The author thanks Fr\'{e} Vercauteren for drawing his attention to the relevance of such results. We don't go into all the details, but the main steps needed for this estimate are the following. Suppose $a=1$, and we only look at the dependency on $n$. As before we compute $r(\G)^{\chi_1}F(\G)$ in time $\Ot(n^2)$, and some Teichm\"{u}ller modulus $\psi(z)$. Let $\og_1, \ldots, \og_k$ be the parameters for which we want to calculate the zeta function. Computing all the Teichm\"{u}ller lifts $\g_1,\ldots,\g_k$ takes $\Ot(n^3)$ time. For computing the matrices $\mathcal{F}_{\g_i}$ we need $F(\g_i^{\sigma^{2^t}})$ for $t=0\ldots \lfloor\log_2 n\rfloor$, hence if we can find all the $\alpha(\g_i^{\sigma^\ell})$ for some $\ell=\O(n)$ and an analytically continuated element $\alpha$ of $F(\G)$ in time $\Ot(n^3)$ we are done.

This is where fast multipoint evaluation pops up. Indeed, computing $\g_i^{\sigma^{\ell}}$ again requires only $\Ot(n^2)$ for each $i$, and the simultaneous substitution of all these values in $\alpha$ takes time $\Ot(n^3)$, which follows from corollary 10.8 in \cite{ModernCompAlg}. The estimate on the memory is clear, as it will certainly not exceed the time requirements.

Note that this result is also applicable to the situation in \cite{HubrechtsHECOdd}, hence for hyperelliptic curves in odd characteristic.

\subsection{Quadratic counting with GNB.}\label{ssec:GNB}

If we work over fields $\F_{q^n}$ where a Gaussian normal basis of type $t$ with $t$ small exists (see e.g. \cite{HandboekHECC}, section 2.3.3.b, and for the existence of such bases \cite{KimParkEA}), then we can make our algorithm quadratic for some well-chosen parameters. Here is an outline of how this works for $t=1$ and $a=1$, which means we have a representation
\[\F_{2^n}\cong\frac{\F_2[x]}{x^n+x^{n-1}+\cdots +x+1}.\]
The same minimal polynomial $(x^{n+1}-1)/(x-1)$ can be used over $\Q_2$ to represent $\Q_{2^n}$, and it is clear that it is a Teichm\"{u}ller modulus. Remark that $x^{n+1}=1$, which makes computing a lot easier. Suppose now that our parameter $\g$ equals some power of $x$, say $x^k$. Note that this is a very strong condition, for there exist only $n+1$ such parameters $\g$. As explained earlier the crucial step is computing $\alpha(\g)^{\sigma^\ell}$ for $\ell=\O(n)$ and $\alpha$ some polynomial of degree $\O(n)$ over $\Q_2$ modulo $2^{\O(n)}$. Now if $\alpha(\G)=\sum_{i=0}^{m} a_i\G^i$, then we have (using a \emph{redundant representation}, i.e. a non-unique form using the generating set $1,x,\ldots,x^n$)
\[\alpha(\g)^{\sigma^\ell}=\alpha(x^{kp^\ell})=\sum_{i=0}^ma_i x^{kip^\ell\bmod n+1},\]
and this last expression is easily evaluated. In conclusion this GNB allows us to compute the zeta function for certain parameters in time $\Ot(n^2)$. Here too we can draw the same conclusions for the odd characteristic case.

\bibliographystyle{amsplain}
\bibliography{bibliography}

\end{document}